\documentclass[12pt]{iopart}
\usepackage[english]{babel}
\usepackage[T1]{fontenc}
\usepackage[latin1]{inputenc}
\usepackage[babel=true]{csquotes} 
\usepackage{amssymb}
\usepackage{amscd}
\usepackage{amsthm}
\usepackage{mathrsfs}
\usepackage[pdftex]{graphicx}
\usepackage{graphicx}

\newcommand{\N}[0]{\ensuremath{\mathbb{N}}}
\newcommand{\E}[0]{\ensuremath{\mathbf{E}}}

\newcommand{\Pe}[0]{\ensuremath{\mathbf{P}}}
\newcommand{\R}[0]{\ensuremath{\mathbb{R}}}
\newcommand{\Z}[0]{\ensuremath{\mathbb{Z}}}

\newtheorem{theo}{\sc{Theorem}}[section]
\newtheorem{lemm}{\sc{Lemma}}[section]
\newtheorem{defi}{\sc{Definition}}[section]

\newtheorem{cor}{\sc{Corollary}}[section]
\newtheorem{rmq}{\sc{Remark}}[section]

\newtheorem{assu}{\sc{Assumption}}[section]
\newtheorem{prop}{\sc{Property}}[section]
\newtheorem*{mainlemm}{\sc{Main Lemma}}

\begin{document}

\title{On hyperbolicity of minimizers for 1D random Lagrangian systems}

\author{Alexandre Boritchev \protect\footnote{Centre de Mathematiques Laurent Schwartz,
Ecole Polytechnique, Route de Saclay
91128 Palaiseau Cedex, France. \\ E-mail: boritchev@math.polytechnique.fr},\ Konstantin Khanin \protect\footnote{Department of Mathematics, University of Toronto, Toronto, Canada. \\ E-mail: khanin@math.toronto.edu}}

%

\begin{abstract}
We prove hyperbolicity of global minimizers for random Lagrangian systems in dimension 1. The proof considerably simplifies a related result in \cite{EKMS}. The conditions for hyperbolicity are almost optimal: they are essentially the same as conditions for uniqueness of a global minimizer in \cite{IK}.
\end{abstract}


\submitto{\NL}

\section{Introduction}
A large body of work on the random forced Burgers equation and Burgers turbulence
in the last 10 years (see \cite{BK} and further references therein) has motivated closely related studies of random
Lagrangian systems \cite{EKMS,IK}. The main object of analysis is a Lagrangian system which depends
smoothly on position $x$ and velocity $v$, but quite irregularly on time $t$:
\begin{equation}
\label{rl}
L^{\omega}(x,v,t)=L_0(x,v) + F^{\omega}(x,t),
\end{equation}
where $F^{\omega}(x,t)$ is a stationary random process in $t$.
The Lagrangian is defined on the tangent bundle $TM$ to a connected $d$-dimensional Riemannian
manifold $M$.
Most rigorous results available at the moment require that $M$ be compact, which will also be the standing assumption in this paper. Since the potential $F^{\omega}(x,t)$ is smooth in $x$ the most natural
continuous time model is given by
\begin{equation}
\label{rp1}
F^{\omega}(x,t)=\sum_{k=1}^K \dot{W}^{\omega}_k(t) F^k(x),
\end{equation}
where $F^k(x)$ are smooth non-random potentials on $M$, and $\dot{W}^{\omega}_k(t)$ are independent white noises. One can also consider \enquote{kicked} models:
\begin{equation}
\label{rp2}
F^{\omega}(x,t)=\sum_{j=-\infty}^{+\infty} F^\omega(j)(x)\delta(t-j),
\end{equation}
where $\lbrace F^{\omega}(j)(x),\ j \in \Z \rbrace$ is a stationary sequence of random potentials. 
We shall assume that potentials $F^{\omega}(j)$ are picked independently for different
$j \in \Z$ according to a given probability distribution $\mu$ on $C^n(M)$, where $n$ is big enough. The Lagrangian
dynamics corresponding to (\ref{rp2}) can be described as follows. For non-integer times $t$ the system evolves
according to a non-random Lagrangian $L_0$, and at integer times $t=j\in \Z$ the velocity 
changes discontinuously:
\begin{eqnarray} \nonumber
&v(j+0)=v(j-0)+ \nabla F^\omega(j)(x).
\end{eqnarray}
Although the two models
(\ref{rp1}) and (\ref{rp2}) look rather different, the theory and results for both cases are parallel.

\

Lagrangian systems (\ref{rl}) are related to random forced Hamilton-
\\
Jacobi equations.
One has to first define the Hamiltonian
\begin{eqnarray} \nonumber
&H^{\omega}(x,p,t)=\max_v{[p\cdot v - L^\omega(x,v,t)]} = H_0(x,p) - F^{\omega}(x,t),
\end{eqnarray}
and then to consider the corresponding Hamilton-Jacobi equation
\begin{equation}
\label{hj}
\phi_t + H^{\omega}(x,\nabla \phi,t)=0.
\end{equation}
One of the most studied cases corresponds to $L_0=v^2/2$. In this case $H_0=p^2/2$ and the Hamilton-Jacobi equation (\ref{hj}) takes the form
\begin{eqnarray} \nonumber
&\phi_t(x,t) + \frac{1}{2}|\nabla \phi|^2- F^{\omega}(x,t)=0.
\end{eqnarray}
Then for the velocity field $v(x,t)=\nabla \phi(x,t)$ one gets the inviscid Burgers equation:
\begin{eqnarray} \nonumber
&v_t(x,t) + (v\cdot \nabla) v(x,t)- \nabla F^{\omega}(x,t)=0.
\end{eqnarray}
Although all the results of this paper hold for any Lagrangian $L_0$ which is convex in $v$ and grows super-linearly as $|v| \to \infty$, below we only consider the case $L_0=v^2/2$.

It is well-known that minimizers for the Lagrangian $L^\omega$ generate the viscosity
solution of the Hamilton-Jacobi equation (\ref{hj}). This connection is especially
useful and important for the study of global solutions, that is solutions for $t \in (-\infty,\infty)$. In order to discuss a global solution one has to fix the value
of the first integral
\begin{equation}
\label{fi}
b=\int_M\nabla \phi (x,t)dx .
\end{equation}
The theory developed in \cite{IK} states that under extremely mild conditions, with probability 1, for every value of the first integral $b\in \R^d$, there exists a unique (up to an additive constant) global solution to the Hamilton-Jacobi equation.  This unique global solution can be viewed as a stationary solution. It plays the role of a global attractor for the dynamics corresponding to the Cauchy problem for the Hamilton-Jacobi equation. Under additional assumptions of non-degeneracy one can also
prove that for every value of $b\in \R^d$, with probability 1, there exists a unique global minimizer for the Lagrangian $L^\omega$ (see \cite {IK}). A global minimizer can be defined
as a smooth curve $\gamma: (-\infty, \infty) \to M$ such that for any compact perturbation $\tilde\gamma$ the difference between Lagrangian actions corresponding to $\tilde\gamma$ and to $\gamma$ is non-negative. Namely, if $\tilde\gamma-\gamma$ is supported on $[T_1,T_2]$, then
\begin{eqnarray} \nonumber
&A^{\omega,b}(\tilde\gamma) - A^{\omega,b}(\gamma)= \int_{T_1}^{T_2}L^\omega(\tilde\gamma,\dot{\tilde\gamma}-b,t)dt 
-\int_{T_1}^{T_2}L^\omega(\gamma,\dot{\gamma}-b,t)dt \geq 0.
\end{eqnarray}
It is expected that
the global minimizer is a hyperbolic trajectory of the Lagrangian flow. Unfortunately such a result is not available at present in the multi-dimensional case $d>1$. In our view hyperbolicity of the global minimizer is one of the most important open problems
in the theory of random Lagrangian systems on compact manifolds. 
In the one-dimensional case hyperbolicity was established in \cite{EKMS}. However the proof
in \cite{EKMS} is unnecessarily complicated and conditions are too restrictive. 
In this paper we present a new proof which is both elementary and conceptual.
Here, conditions for hyperbolicity are almost the same as the conditions
for uniqueness of a global minimizer (see \cite{IK}). This is another important advantage of the approach used in this article. 

The following property is crucial for establishing hyperbolicity of the global
minimizer. 
Define first backward minimizers as minimizers on semi-infinite time intervals $(-\infty, t]$ with one end point at $t$ fixed. They can be defined in the same way as global minimizers. Now consider all backward minimizers which originate at time $t$,
and denote by $\Omega_{s,t}$ the set of all points $x$ which are reached by some backward minimizer at time $s \leq t$. We prove that the diameter of $\Omega_{s,t}$ tends to zero exponentially as $t \to \infty$. This property implies hyperbolicity by the standard
argument, which also allows to construct corresponding stable and unstable manifolds. We shall not discuss these issues in the present paper and refer the readers to \cite {EKMS}. Instead, here we shall only deal with the key shrinking property
formulated above.

\

We finish this section with several general remarks. First, we want to emphasize 
the importance of hyperbolicity of the global minimizer. It immediately
implies many fundamental properties of the global solution to
the Hamilton-Jacobi equation, such as piecewise smoothness, 
exponential rate of convergence to the global solution, and many others. It also allows to study the structure of singularities (shocks) (see \cite {BK}).  

\

Our second remark is related to a general problem of hyperbolicity of minimizers
for generic non-random Lagrangian systems. This is one of the central problems of the
Aubry-Mather theory. Randomness is another way to introduce the notion of genericity. In this setting generic stands for properties which
hold for almost all systems (with probability 1). Note however that in many respects, random
and nonrandom (autonomous, or depending on time periodically) Lagrangian systems are 
very different. In particular, all number-theoretical
aspects of the Aubry-Mather theory disappear in the random case.

\

Finally, we want to say a few words about the non-compact case. At present there are almost no
rigorous results in that setting. It is believed that if the system exhibits any
form of translation invariance, global minimizers do not exist. However, it is likely
that backward minimizers do exist, and the study of their asymptotic scaling properties is an  extremely interesting and important problem. 


\section{Hyperbolicity assumptions and main results}

We begin by formulating the assumptions on potentials:

\begin{assu} \label{C1}
In the \enquote{kicked} case, we assume the following.
\\
(i)\ The kicks at integer times $j$ are of the form
$$
F^{\omega}(j)=\sum_{k=1}^{K}{c_k^{\omega}(j) F^k},
$$
where $F^k$ are smooth potentials on $S^1=\R/\Z$. The random vectors
\\
$(c_k^{\omega}(j))_{1 \leq k \leq K}$ are independent identically distributed $\R^K$-valued random variables. Their distribution on $\R^K$, denoted by $\mu$, is assumed to be absolutely continuous with respect to the Lebesgue measure.
\\
(ii)\ 
$0$ belongs to $Supp\ \mu$.
\\
(iii)\ The mapping from $S^1$ to $\R^K$ defined by
$$
x \mapsto (F^1(x),...,F^K(x))
$$
is an embedding.
\end{assu}

\begin{rmq}
Let $g$ be the function defined by
$$
g(c_1,...,c_K)=\sum_{k=1}^{K}{c_k F^k}.
$$
We denote by $\nu$ the corresponding push-forward measure
$$
\nu=g_{*}(\mu)
$$
on a smooth Sobolev space. The assumption $0 \in Supp\ \mu$ can then be replaced by the slightly weaker assumption $0 \in Supp\ \nu$.
\end{rmq}

\begin{assu} \label{C2}
In the case of the white force potential, we assume the following.
\\
(i)\ The forcing has the form
$$
F^{\omega}(x,t)=\sum_{k=1}^{K}{\dot{W}_k^{\omega}(t) F^k(x)},
$$
where $F^k$ are smooth potentials on $S^1$, and $\dot{W}^{\omega}_k$ are independent white noises, i.e. weak time derivatives of independent Wiener processes $W_k^{\omega}(t)$.
\\
(ii)\ The mapping from $S^1$ to $\R^K$ defined by
$$
x \mapsto (F^1(x),...,F^K(x))
$$
is an embedding.
\end{assu}
\smallskip
\indent
We denote by $G$ an antiderivative in time of the forcing:
$$
G^{\omega}(x,t)=\sum_{k=1}^{K}{W_k^{\omega}(t) F^k(x)},
$$
where $W_k^{\omega}(t)$ are independent standard Wiener processes with $W_k^{\omega}(0)=0$. Since we will only consider time differences of $G$, the particular choice of antiderivative has no importance.
\\ \indent
In both cases, $F^{\omega}$ will be abbreviated as $F$, and in the white force case $F(\cdot,t)$ will be abbreviated as $F(t)$, and similarly for $G$.

\begin{rmq}
The embedding conditions are consistent with the condition for uniqueness of the global minimizer (see \cite{IK}). In the \enquote{kicked} case, the condition for uniqueness in \cite{IK} is slightly weaker: the map $x \mapsto (F^1(x),...,F^K(x))$ is  only required to be one-to-one. However, we need to assume the embedding to prove hyperbolicity.
\end{rmq}

The following property, called the \textit{separation property}, plays a crucial role in our construction.

\begin{prop}
There exist $\alpha_0>0$, three pairwise disjoint open intervals $J_i$, $i=1,2,3$, and three potentials $\tilde{F}_i$, $i=1,2,3$
with the following properties.
\\
1) In the \enquote{kicked} case, we have $\tilde{F}_i \in Supp\ \nu$ for every $i$. In the white force case, each $\tilde{F}_i$ is a linear combination of the $F^k$.
\\
2)\ Each of the functions $-\tilde{F}_i$ reaches its minimum, denoted by $m_i$, at a single point $x_i$.
\\
3)\ For every $\alpha,\ 0<\alpha \leq \alpha_0$, there exist three open intervals $I_i(\alpha)$, $I_i \subset J_i,\ i=1,2,3$ such that
\begin{eqnarray} \nonumber
&\tilde{F}_i(S^1-I_i) \subset (-\infty,-m_i-\alpha].
\end{eqnarray}
\end{prop}

Note that for every $i$ and $\alpha$, the point $x_i$ where $\min (-\tilde{F}_i)$ is reached belongs to $I_i$.

\begin{lemm} \label{implysep}
Assumptions \ref{C1} or \ref{C2} imply the separation property.
\end{lemm}

\textbf{Proof of Lemma \ref{implysep}:}
\\
\textbf{\enquote{Kicked} case:} We start by showing that, for Lebesgue-a.e. vector $(c_j)_{1 \leq j \leq K}$, the maximum of
$$
\sum_{j=1}^{K}{c_j F^j(x)}
$$
is reached at a single point $x \in S^1$. This follows from a rather standard argument (see \cite[Corollary 5]{IK}). Indeed, the function
$$
\Phi: (c_1,\dots,c_K) \mapsto \max_{x \in S^1} \sum_{j=1}^{K}{c_j F^j(x)}
$$
is Lipschitz and therefore differentiable a.e., with respect to the
\\
Lebesgue measure $\mu_{Leb}$. On the other hand, at a point of differentiability of $\Phi$,
$$
\nabla \Phi(x_{max}) =(F^1(x_{max}),\dots,F^K(x_{max}))
$$
for every point of maximum $x_{max}$. Hence the embedding assumption \ref{C1} (iii) implies that the point of maximum is unique. Since $\mu$ is absolutely continuous with respect to $\mu_{Leb}$, the maximum uniqueness set $O_1 \subset \R^K$ has full $\mu$-measure.
\\ \indent
Furthermore, by the Lebesgue points theorem \cite[Theorem 7.7]{Rud}, $c=(c_j)_j$ is a Lebesgue point for the density
$$
q=\frac{d \mu}{d \mu_{Leb}}
$$
on a set $O'$ of full $\mu_{Leb}$-measure, and thus of full $\mu$-measure.
\\ \indent
Denote by $O_2 \subset O'$ the set of Lebesgue points $c$ for $q$ such that $q(c)>0$. By definition, they belong to $Supp\ \mu$, and $O_2$ has full $\mu$-measure.
\\ \indent
Now consider $c^1=(c^1_j)_j \in O_1 \cap O_2$.
Denote by $x_1$ the point where the maximum of $\tilde{F}_1= \sum_{j=1}^{K}{c^1_j F^j}(x)$ is reached: $x_1=argmax\ \tilde{F}_1$.
\\ \indent
Denote by $V$ the set of vectors $(c_j)_j$ such that
$$
\sum_{j=1}^{K}{c_j \frac{d F^{j}}{dx}(x_1)} \neq 0.
$$
Denote by $B_n$ the open ball with radius $1/n$ centered at $c^1$. We will also need $B^{'}_n=B_n \cap (c^1+V) \cap O_1 \cap O_2$. 
By the embedding assumption \ref{C1} (iii), $B_n \cap (c^1+V)$ is just $B_n$ itself with a removed hyperplane. Thus, since $\mu$ is continuous with respect to $\mu_{Leb}$, we have $\mu(B^{'}_n)= \mu(B_n).$
\\ \indent
Using \cite[Theorem 7.7]{Rud} one more time, we obtain that there exists a constant $N_0$ such that for $n \geq N_0$,
\begin{eqnarray} \nonumber
&\mu(B^{'}_n) =\mu (B_n) \geq \frac{q(c^1)}{2} \mu_{Leb} (B_n) > 0.
\end{eqnarray}
\indent
On the other hand, for small enough $\epsilon>0$ there exists $N_1(\epsilon)$ such that for $n \geq N_1$, if $(c_j)_j \in B^{'}_n$, then $\sum_{j=1}^{K}{c_j F^j}$ reaches its (unique) maximum in a point of the $\epsilon$-neighbourhood of $x$ different from $x$ itself. Considering a smaller neighbourhood at each step, this argument can be repeated any finite number of times. It enables us to construct any number of potentials contained in $Supp\ \nu$ and attaining their respective maxima at different points: three suffice for our purposes. Denote them by $\tilde{F_1}, \tilde{F_2}, \tilde{F_3}$. Let $J_1, J_2, J_3$ be three non-intersecting open intervals around their respective points of maximum. Take as $\alpha_0$ the minimum of $\max(\tilde{F_i})-\max(\tilde{F_i}|_{S^1-J_i})$. It is obvious that for any $\alpha \in (0,\alpha_0]$ we can construct the required intervals $I_i(\alpha)$.
\\ \indent
\textbf{White force case:} The proof follows the same lines, but is much simpler since measure-theoretic arguments are trivialised. $\square$

\begin{defi}
Consider a closed subset $Z$ of $S^1$. Let $m(Z)$ denote the maximal length of a connected component of $S^1-Z$. We define the diameter of $Z$ as
$$
d(Z)=1-m(Z).
$$ 
\end{defi}

\indent
The diameter of $Z$ can be thought of as the minimal length of an interval on $S^1$ containing $Z$.
\smallskip
\\ \indent
In what follows we use the function $\psi^{\omega}$, either deterministic or random, as an initial condition at time $s$. Everywhere below, the value of the first integral $b$ (see (\ref{fi})) is fixed. For simplicity, we do not indicate dependence on $b$ in our notation.

\begin{defi}
For a given value of $b \in \R$, a curve $\gamma_{s,t}^{y,x}(\tau)$ is a minimizer if it minimizes the action
\begin{eqnarray} \nonumber
&A(\gamma)=\frac{1}{2} \int\limits_{s}^{t}{(\dot{\gamma}(\tau)-b)^2 d \tau} + \sum_{n \in [s,t)} {\Big(-F(n)(\gamma(n))\Big)}
\end{eqnarray}
in the \enquote{kicked} case and the action
\begin{eqnarray} \nonumber
A(\gamma) =&\frac{1}{2} \int\limits_{s}^{t}{(\dot{\gamma}(\tau)-b)^2 d \tau}
\\ \nonumber
&+\int\limits_{s}^{t} { \Bigg( \dot{\gamma}(\tau) \Big(\frac{\partial G}{\partial x}(\gamma(\tau),\tau)-\frac{\partial G}{\partial x}(\gamma(\tau),t) \Big) \Bigg) d \tau }
\\ \nonumber
&+ \Big( G(\gamma(s),s)-G(\gamma(s),t) \Big)
\end{eqnarray}
in the white force case, respectively, over all absolutely continuous
\\
curves with endpoints $x$ at time $t$ and $y$ at time $s$.
\end{defi}

\begin{defi}
For any time interval $[s,t]$ and any continuous function $\psi: S^1 \rightarrow \R$, a curve $\gamma_{s,t,\psi}^{x}(\tau): [s,t] \mapsto S^1$ is a $\psi$-minimizer if it minimizes 
$A(\gamma)+\psi(\gamma(s))$
over all absolutely continuous curves with  endpoint $x$ at time $t$.
\end{defi}

\begin{defi}
For $-\infty<r<s \leq t<+\infty$ and for a fixed function $\psi(\cdot,r): S^1 \rightarrow \R$, let $\Omega_{r,s,t,\psi}$ be the set of points reached, at the time $s$, by $\psi$-minimizers on $[r,t]$:
\begin{eqnarray} \nonumber
&\Omega_{r,s,t,\psi}=\lbrace \gamma_{r,t,\psi}^x(s),\ x \in S^1\rbrace.
\end{eqnarray}
\end{defi}

\begin{rmq}
In what follows, the initial condition $\psi$ will always be fixed, while $t$ will increase to $+\infty$. It is important that we shall consider both deterministic and random initial conditions $\psi$. In the latter case, $\psi$ should be measurable with respect to the past $\sigma$-algebra $\mathcal{B}_r=\mathcal{B}_{(-\infty,r]}$, which is defined in a standard way. It is important  to take $r$ smaller then $s$. Everywhere below, we set $r=s-1$. To simplify notation, $\Omega_{s-1,s,t,\psi}$ will be denoted by $\Omega_{s,t}$.
\end{rmq}

It is well-known that $\Omega_{s,t}$ is a closed set. Obviously, $\Omega_{s,t_1} \supseteq \Omega_{s,t_2}$ for all $s \leq t_1 \leq t_2$. It follows that $t \mapsto d(\Omega_{s,t})$ is a non-increasing function.
\smallskip \\ \indent
We are now able to formulate the main results of this paper which are the following theorem and its corollary. Both results hold for a given value of $b \in \R$. However, all constants are uniformly bounded if $b$ stays bounded. It is easy to see that in the \enquote{kicked} case, $b$ is effectively defined modulo $1$, since the action is invariant under the transformation $(b,\gamma) \mapsto (b+1,\gamma+t)$. Thus in this case all constants are uniformly bounded for all $b$.

\begin{theo} \label{main}
Assume that the separation property holds. Then there exist constants $\lambda,\tilde{C}>0$ such that if $-\infty < s \leq t < +\infty$, then
\begin{eqnarray} \nonumber
&\E(d(\Omega_{s,t})) \leq \tilde{C} \exp(-\lambda(t-s)),
\end{eqnarray}
where $\E(\cdot)$ stands for the expectation with respect to the distribution of potentials.
\end{theo}

\begin{cor} \label{maincor}
Assume that the separation property holds. Fix 
\\
$s \in \R$. Then, for a.e. $\omega$, there exists a random constant $\tilde{C}(s,\omega)>0$ such that
\begin{eqnarray} \nonumber
&d(\Omega_{s,t}) \leq \tilde{C}(s,\omega) \exp(-\lambda(t-s)/2),\quad t \geq s.
\end{eqnarray}
Here, $\lambda$ is the same as in Theorem \ref{main}.
\end{cor}

As we have already pointed out in the introduction, Corollary \ref{maincor} implies hyperbolicity (see \cite{EKMS} for details). The following lemma, called the \textit{main lemma}, is proved in Section \ref{mainsection}: the proof is quite involved, with additional technical difficulties in the white force case.

\begin{mainlemm}
Assume that the separation property holds. Fix $b \in \R$. Then there exist constants $c,T>0$ such that if $-\infty < s \leq t < +\infty$, then the following inequality holds a.s.:
\begin{eqnarray} \nonumber
&\Pe \Bigg(d(\Omega_{s,t+T}) \leq \frac{d(\Omega_{s,t})}{2}\ |\ \mathcal{B}_t \Bigg) \geq c.
\end{eqnarray}
\end{mainlemm}

We finish this section by deriving Theorem \ref{main} and Corollary \ref{maincor} from the main lemma.
\\ \bigskip
\\
\textbf{Proof of Theorem \ref{main}
:}
Consider the function
\begin{eqnarray} \nonumber
&d(t)=\E(d(\Omega_{s,t})) \exp(\lambda(t-s)),
\end{eqnarray}
where $\lambda$ is a fixed positive number, chosen later.
\\ \indent
Since $t \mapsto d(\Omega_{s,t})$ is non-increasing, the main lemma implies that
\begin{eqnarray} \nonumber
&\E(d(\Omega_{s,t+T})) \leq c\ \frac{\E(d(\Omega_{s,t}))}{2}  + (1-c) \E(d(\Omega_{s,t})).
\end{eqnarray}
Thus
\begin{eqnarray} \nonumber
&d(t+T) \leq \exp(\lambda T)\ \Big( 1-\frac{c}{2} \Big) d(t).
\end{eqnarray}
Now put
\begin{eqnarray} \nonumber
&\lambda = -\frac{1}{T} \ln \Big( 1-\frac{c}{2} \Big) .
\end{eqnarray}
It follows that $d(t+T) \leq d(t)$. But  $d(s) = 1$. Therefore, for $t \in s+T \N$, we have $d(t) \leq 1$.
Consequently, since $t \mapsto d(\Omega_{s,t})$ is non-increasing, we have
\begin{eqnarray} \nonumber
&\E(d(\Omega_{s,t})) \leq \tilde{C} \exp(-\lambda(t-s)),\quad t \geq s,
\end{eqnarray}
with $\tilde{C}=\exp(\lambda T)=(1-\frac{c}{2})^{-1}$. This proves the theorem's assertion. $\square$
\medskip

\textbf{Proof of Corollary \ref{maincor} assuming Theorem \ref{main}:}
In the same way as in the previous proof, it is enough to prove the statement for $t \in s+ T \N$.
By Theorem $\ref{main}$ and Chebyshev's inequality, for every $X>0$,
\begin{eqnarray} \nonumber
&\Pe(d(\Omega_{s,s+nT}) \geq X \exp(-\lambda n T/2)) \leq \frac{\tilde{C}}{X} \exp(-n\lambda T/2),\quad n \geq 0.
\end{eqnarray}
An application of the Borel-Cantelli lemma ends the proof. $\square$
\bigskip 
\\

\section{Proof of the main lemma}\label{mainsection}

For all $s<t$, let us define a map $S_s^t$ from $S^1$ to $S^1$, which can be viewed as a coordinate projection at time $t$ of the generalized Lagrangian flow corresponding to the Burgers equation. It certainly depends on the initial condition $\psi$ at time $s-1$.
\\ \indent
If, at time $s$, a point $y$  belonging to $S^1$ is reached by a $\psi$-minimizer on $[s-1,t]$ starting in $x$ at time $t$, then $S_s^t(y)$ is equal to the point $x$. Note that such an $x$ is unique, since minimizers on the time interval $[s-1, t]$ cannot intersect outside of endpoints $s-1$ and $t$. 
\\ \indent
If a point $y$ is not reached by such a $\psi$-minimizer, then it belongs to a closed interval corresponding to a shock at time $t$. In this case $S_s^t(y)$ is equal to the corresponding shock position. To define an interval at time $s$ corresponding to a shock at time $t$, one has to consider rightmost and leftmost minimizers originating at $(x,t)$. Intersections of those minimizers with $S^1 \times \lbrace s \rbrace$ generate a space interval of points absorbed by the shock $(x,t)$. It is easy to see that every point $(y,s)$ is reached by a minimizer or belongs to a shock interval generated by a uniquely defined shock.
\\ \indent
Note that some points may correspond to both cases considered above. Namely, points corresponding to minimizers which originate from the shock positions. However, even in this case the map $S_s^t$ is still uniquely defined.

\subsection{Proof in the \enquote{kicked} case}

Put
\begin{equation} \label{C}
C=3\ \Big( \max_{i \in {1,2,3}}{\Vert\tilde{F}_i\Vert_{C^1}} + 1 \Big).
\end{equation}
Then put
\begin{equation} \label{alpha}
\alpha=\min \Bigg( \alpha_0,\frac{1}{10C} \Bigg)
\end{equation}
(see the separation property for the definition of $\alpha_0$.) We keep in mind that $\alpha < 1/30$.
\\ \indent
Consider integers
\begin{equation} \label{NN}
N' \in \Big( 2+\frac{1}{\alpha^3},\ \frac{2}{\alpha^3} \Big);\ N \in \Big( \frac{2}{\alpha^{10}}+1,\ \frac{4}{\alpha^{10}} \Big).
\end{equation}
Denote by $E_1$ the event
\begin{equation} \label{0Nkick}
\Vert F(t+k)\Vert_{\infty} \leq \alpha^{20},\quad 0 \leq k \leq N-1.
\end{equation}
By Assumption \ref{C1} $(ii)$ the zero potential belongs to $Supp\ \nu$. It follows that $E_1$ has positive probability.
\\ \indent
Put $l=1-d(\Omega_{s,t})$. If $\Omega_{s,t} \neq S^1$, consider a connected component $(y_1,y_2)$ of $S^1-\Omega_{s,t}$ which has maximal length $l$. Let $y_3$ be the center of $(y_1,y_2)$, and let $y_4$ be the point diametrically opposite to $y_3$. If $\Omega_{s,t}=S^1$, let $y_3$ and $y_4$ be any pair of diametrically opposite points in $S^1$. Then consider $z_1=S_s^{(t+N)}y_3$ and $z_2=S_s^{(t+N)}y_4$. Since the $J_i$ (see the separation property for their definition) are pairwise disjoint, one of the $J_i$ has an empty intersection with one of $[z_1,z_2]$ and $[z_2,z_1]$.
Without loss of generality, we may suppose that $[z_1,z_2] \cap J_1=\varnothing$.
\\ \indent
Now consider the straight line defined by
\begin{eqnarray} \nonumber
&\gamma(\tau)=x+b(\tau-t-N),\quad \tau \in [t+N,t+2N-1]
\end{eqnarray}
for some $x \in S^1$.
\\ \indent
We claim that there exist (at least) $N'$ different integers $0=n_0 < \dots < n_{N'-1} \leq N'-1$ such that we have
\begin{equation} \label{njaux}
\max_{j,j' \in [0,N-1]} |\gamma(t+N+n_j)-\gamma(t+N+n_{j'})| \leq \alpha^7.
\end{equation}
Indeed, by the pigeonhole principle, since $N' \leq \alpha^7 N$, there exist integers $0 \leq \tilde{n}_0 < \dots < \tilde{n}_{N'-1} \leq N-1$ such that
\begin{eqnarray} \nonumber
&\max_{j,j' \in [0,N'-1]} |\gamma(t+N+\tilde{n}_j)-\gamma(t+N+\tilde{n}_{j'})| \leq \alpha^7.
\end{eqnarray}
Then it suffices to take, for every $j$, $n_j=\tilde{n}_j-\tilde{n}_0$.
\\ \indent
By definition of $C$ and $\alpha$, (\ref{njaux}) yields that
\begin{eqnarray} \nonumber
&\max_{j,j' \in [0,N'-1]} |\tilde{F}_1(\gamma(t+N+n_j)) - \tilde{F}_1(\gamma(t+N+n_{j'})) |
\\ \indent \label{nj}
&\leq \alpha^7 \Vert\tilde{F}_1\Vert_{C^1} \leq \alpha^6/10.
\end{eqnarray}
Now consider the event $E_2$ defined by the system of inequalities:
\begin{equation} \label{N2Nkick}
\cases{
\Vert F(t+N+n_j)-\tilde{F}_1\Vert_{\infty} \leq \alpha^{20},\quad 0 \leq j \leq N'-1.
\\
\Vert F(t+N+k)\Vert_{\infty} \leq \alpha^{20},
\\
k \in [0,N-1]-\lbrace n_0,\dots,n_{N'-1} \rbrace.
}
\end{equation}
Since $\tilde{F}_1$ and $0$ belong to $Supp\ \nu$, this event (independent from $E_1$) also has positive probability.
\\ \indent
It remains to prove that for $\omega \in E_1 \cap E_2$ all minimizers on $[t, t+2N]$ pass through $I_1(\alpha)$ at time $t+N$, which follows from Lemma \ref{kick1} and Lemma \ref{kick2}. Indeed, 
if this statement holds, no such minimizers can pass through $[z_1,z_2]$ at $t+N$, since $I_1(\alpha) \subset J_1$.  Consequently all the points that are in $[y_3,y_4]$ at time $s$ will not be reached by minimizers originating at time $t+2N$.  In particular, it follows that $[y_3,y_4]$ is contained in an interval generated by some shock at time $t+2N$. Therefore $(y_1,y_2) \cup [y_3,y_4] = (y_1,y_4]$ is contained in a connected component of $S^1-\Omega_{s,t+2N}$.
Thus
\begin{eqnarray} \nonumber
&d(\Omega_{s,t+2N}) \leq 1-\frac{1+l}{2} = \frac{1}{2} d(\Omega_{s,t})
\end{eqnarray}
with a positive conditional probability which equals at least $\Pe(E_1)\Pe(E_2)$. This proves the lemma's assertion. $\square$

\begin{lemm} \label{kick1}
Assume that $\omega \in E_2$. Then for every minimizer $\gamma_1$ on $[t+N,t+2N]$ there exists $j,\ 1 \leq j \leq N'-1$,  such that
\begin{eqnarray} \nonumber
&-\tilde{F}_1(\gamma_1(t+N+n_j)) \leq m_1+\alpha^2.
\end{eqnarray}
\end{lemm}

\textbf{Proof:} We argue by contradiction. Suppose that there exists a minimizer $\gamma_1$ on $[t+N,t+2N]$ such that
\begin{equation} \label{kick1contra}
-\tilde{F}_1(\gamma_1(t+N+n_j)) > m_1+\alpha^2,\quad 1 \leq j \leq N'-1.
\end{equation}
Consider a curve $\gamma_2$ with the same endpoints as $\gamma_1$, linear on intervals $[t+N+k,t+N+k+1]$.
Moreover we suppose that $\gamma_2 = x_1+b(\tau-t-N)$ on $[t+N+n_1,t+N+n_{N'-1}]$ ($x_1$ being the point where $\tilde{F}_1$ reaches its maximum), and that $|\dot{\gamma}_2-b| \leq 1/2n_1$ and $|\dot{\gamma}_2-b| \leq 1/2(N-n_{N'-1})$ on the extremal intervals $[t+N,t+N+n_1]$ and $[t+N+n_{N'-1},t+2N]$, respectively.
\\ \indent
From now on, for a curve $\gamma$ we denote $\dot{\gamma}-b$ by $ \dot{\gamma}^b$. We recall that the \enquote{kicked} case action $A$ for $\gamma: [t_1,t_2] \rightarrow S^1$ equals
\begin{eqnarray} \nonumber
&A(\gamma)=\frac{1}{2} \int\limits_{t_1}^{t_2}{(\dot{\gamma}^b(\tau))^2 d \tau} - \sum_{n \in [t_1,t_2)} {F[\gamma(n)]}.
\end{eqnarray}
The first part of the right-hand side, corresponding to the kinetic energy, will be denoted by $A^k$. The remaining part, corresponding to the potential energy, will be denoted by $A^p$. We observe that
\begin{equation} \label{Chasles}
A^k(\gamma|_{[t_1,t_3]})=A^k(\gamma|_{[t_1,t_2]})+A^k(\gamma|_{[t_2,t_3]}),
\end{equation}
and similarly for $A^p$. We have
\begin{eqnarray} \nonumber
&A^k(\gamma_1) \geq 0;\ A^k(\gamma_2) \leq \frac{1}{4}.
\end{eqnarray}
On the other hand, using the inequalities (\ref{nj}-\ref{kick1contra}), we get
\begin{eqnarray} \nonumber
&A^p(\gamma_1) \geq (N'-1) ( m_1+\alpha^2-\alpha^{20} ) - (N-N') \alpha^{20}- F(\gamma(t+N)),
\\ \nonumber
&A^p(\gamma_2) \leq (N'-1) ( m_1+\alpha^6/10+\alpha^{20} ) + (N-N') \alpha^{20}- F(\gamma(t+N)).
\end{eqnarray}
Therefore, by (\ref{alpha}-\ref{NN}), we get
\begin{eqnarray} \nonumber
A(\gamma_1)-A(\gamma_2) &= A^k(\gamma_1)  -A^k(\gamma_2) + A^p(\gamma_1) - A^p(\gamma_2)
\\ \nonumber
&\geq -\frac{1}{4}+(N'-1)(\alpha^2-\alpha^6/10)-2(N-1) \alpha^{20}
\\ \nonumber
&\geq -\frac{1}{4}+\alpha^{-1}-\frac{\alpha^3}{10}-8\alpha^{10}>0.
\end{eqnarray}
Thus we have a contradiction with the fact that $\gamma_1$ is a minimizer. This proves the lemma's assertion. $\square$

\begin{lemm} \label{kick2}
Assume that $\omega \in E_1 \cap E_2$. For some $j,\ 1 \leq j \leq N'-1$, consider a minimizer $\gamma_1$ on $[t,t+N+n_j]$ such that  $y=\gamma_1(t+N+n_j)$ satisfies:
\begin{eqnarray} \nonumber
&-\tilde{F}_1(y) \leq m_1+\alpha^2.
\end{eqnarray}
Then we have
\begin{eqnarray} \nonumber
&\gamma_1(t+N) \in I_1(\alpha).
\end{eqnarray}
\end{lemm}

\textbf{Proof:} We argue by contradiction, supposing that $\gamma_1(t+N) \notin I_1(\alpha)$.
We may also assume that
\begin{eqnarray} \nonumber
&-\tilde{F}_1(\gamma_1(t+N+n_{j'})) > m_1+\alpha^2,\quad 1 \leq j' < j.
\end{eqnarray}
Indeed, otherwise we could consider a smaller value of $j$.
In the same way as previously, we want to prove that $\gamma_1$ cannot be a minimizer, and we consider a curve $\gamma_2$ with the same endpoints as $\gamma_1$. Namely, we suppose that $\gamma_2$ satisfies $\dot{\gamma}_2^b=0$ between $t+N$ and $t+N+n_j$, $\gamma_2$ is linear between $t$ and $t+N$, and moreover $|\dot{\gamma}_2^b| \leq 1/2N$. We have the inequalities
\begin{eqnarray} \nonumber
&A^k(\gamma_1) \geq 0;\ A^k(\gamma_2) \leq \frac{1}{8N}.
\end{eqnarray}
On the other hand, using the separation property, (\ref{0Nkick}), (\ref{nj}), and (\ref{N2Nkick}), we get
\begin{eqnarray} \nonumber
A^p(\gamma_1) &\geq -N \alpha^{20} + (m_1+\alpha-\alpha^{20}) + (j-1) (m_1+\alpha^2-\alpha^{20})
\\ \nonumber
&-(n_j-j)\alpha^{20},
\\ \nonumber
A^p(\gamma_2) &\leq N \alpha^{20} +j ( m_1+\alpha^2+\alpha^6/10+\alpha^{20} ) +(n_j-j)\alpha^{20}.
\end{eqnarray}
Therefore, by (\ref{alpha}-\ref{NN}), we obtain that
\begin{eqnarray} \nonumber
&A(\gamma_1)-A(\gamma_2) \geq -\frac{1}{8N}+\alpha-\alpha^2 - N' \alpha^6/10 - 4N \alpha^{20} > 0.
\end{eqnarray}
Again, we have a contradiction. This proves the lemma's assertion. $\square$
\medskip

\subsection{Proof in the white force case}

The scheme of the proof is very similar to the one in the \enquote{kicked} case. The major differences are auxiliary lemmas which are technically more involved and the conditions on the forcing, in some way much more restrictive.
\\ \indent
The constants $C,\alpha,N',N$ are the same as in the proof of the \enquote{kicked} case, with the exception that now
\begin{equation} \label{alphabis}
\alpha=\min \Bigg( \alpha_0,\frac{1}{10C},\frac{1}{10(b+1)^2} \Bigg),
\end{equation}
and that the definitions of $N'$ and $N$ change accordingly. Denote by $E_1$ the event
\begin{equation} \label{0Nwhite}
\sup_{t_1,t_2 \in [t,t+N]} \Vert G(t_1)-G(t_2)\Vert_{C^1} \leq \alpha^{40}.
\end{equation}
By classical properties of the Wiener process, $E_1$ has positive probability, uniformly in $t$.
\\  \indent
Now we proceed exactly in the same way as in the \enquote{kicked} case, supposing with the same notation and without loss of generality that $[z_1,z_2] \cap J_1=\varnothing$.
\\  \indent
We assume that for every $j,\ j \in [0,N'-1]$ (we take $n_{N'}=N$), $G$ satisfies: 
\begin{equation} \label{N2Nwhite}
\cases{
\Vert G(t+N+n_{j+1})-G(t+N+n_j)-\tilde{F}_1\Vert_{C^1} \leq   \alpha^{40}. 
\\
\Vert G(t+N+n_{j+1})-G(t+N+n_j+\tau)\Vert_{C^1} \leq   \alpha^{40},
\\
\tau \in [  \alpha^{40}, n_{j+1}-n_j].
\\
\Vert G(t+N+n_{j}+\tau) - G(t+N+n_j+\tau') \Vert_{C^1}
\\
\leq \frac{3}{2} \Vert \tilde{F}_1 \Vert_{C^1} \leq \frac{C}{2},\ \tau,\tau' \in [0, n_{j+1}-n_j].
}
\end{equation}
This event, denoted by $E_2$, has positive probability and is independent from $E_1$.
\\ \indent
Finally, in the same way as in the \enquote{kicked} case, the lemma's assertion follows from Lemma \ref{white1} and Lemma \ref{white2}.
$\square$

\begin{lemm} \label{white1}
Consider a minimizer $\gamma_1$ on $[t+N,t+2N]$. Then, if $\omega \in E_2$, we have
\begin{equation} \label{inwhite1}
-\tilde{F}_1(\gamma_1(t+N+n_j)) \leq m_1+\alpha^2
\end{equation}
for some $j,\ 1 \leq j \leq N'-1$.
\end{lemm}

\textbf{Proof:} As previously, we argue by contradiction, considering a minimizer $\gamma_1$ on $[t+N,t+2N]$ such that (\ref{inwhite1}) does not hold for any $j,\ 1 \leq j \leq N'-1$. We recall that the action is given by:
\begin{eqnarray} \nonumber
A(\gamma) =&\frac{1}{2} \int\limits_{t_1}^{t_2}{ \dot{\gamma}^b(\tau)^2 d \tau}
\\ \nonumber
&+ \int\limits_{t_1}^{t_2}{ \Bigg(  \dot{\gamma}(\tau) \Big(\frac{\partial G}{\partial x}(\gamma(\tau),\tau)-\frac{\partial G}{\partial x}(\gamma(\tau),t_2) \Big) \Bigg) d \tau}
\\  \nonumber
&+ \Big( G(\gamma(t_1),t_1)-G(\gamma(t_1),t_2) \Big).
\end{eqnarray}
The first term of the right-hand side, i.e. the kinetic energy, will be denoted by $A^1$. The second and the third terms, whose sum is the potential energy, will be denoted by $A^2$ and $A^3$, respectively. We observe that $A$ as well as the quantities $A^1$ and $A^{2}+A^{3}$ satisfy a relation of the same type as (\ref{Chasles}). To see it for $A^{2}+A^{3}$, it suffices to write down this sum as a stochastic integral. $A(\gamma|_{[s,t]})$ is denoted by $A_{s,t}(\gamma)$, and similarly for $A^i,\ i=1,2,3$.
\\ \indent
Consider a curve $\gamma_2$ with the same endpoints as $\gamma_1$, defined exactly in the same way as in the proof of Lemma \ref{kick1}. Namely, $\gamma_2 = x_1+b(\tau-t-N)$ on $[t+N+n_1,t+N+n_{N'-1}]$, and $\gamma_2$ is linear on $[t+N,t+N+n_1]$ and on $[t+N+n_{N'-1},t+2N]$ with $|\dot{\gamma}_2^b| \leq 1/2n_1$ and $|\dot{\gamma}_2^b| \leq 1/2(N-n_{N'-1})$, respectively.
\\ \indent
Now, for every $j \in [0,N'-1]$, consider a straight line $\gamma_3$ connecting $\gamma_1(t+N+n_j)$ and $\gamma_1(t+N+n_{j+1})$ 
with constant velocity $\dot{\gamma_3},\ |\dot{\gamma_3}^b| \leq 1/2(n_{j+1}-n_j)$ (we take $n_{N'}=N$). Denote by $R$ the quantity
\begin{eqnarray} \nonumber
&R=A^2_{t+N+n_j,t+N+n_{j+1}}(\gamma_1).
\end{eqnarray}
Since
\begin{eqnarray} \nonumber
&\int{a(\tau)b(\tau) d\tau} \geq - \frac{2\alpha^{20}}{C} \int{\Big( \frac{a(\tau)}{2} \Big)^2 d \tau} - \frac{C}{2\alpha^{20}} \int{b(\tau)^2 d \tau},
\end{eqnarray}
then we have
\begin{eqnarray} \nonumber
R & \geq - \frac{2\alpha^{20}}{C} \int\limits_{t+N+n_j}^{t+N+n_{j+1}}{\Big(\frac{\dot{\gamma_1^b}(\tau)+b}{2}\Big)^2 d \tau} 
\\ \nonumber
&- \frac{C}{2\alpha^{20}} \int\limits_{t+N+n_j}^{t+N+n_{j+1}}{ \Big(\frac{\partial G}{\partial x}(\gamma_1(\tau),\tau)-\frac{\partial G}{\partial x}(\gamma_1(\tau),t+N+n_{j+1}) \Big)^2 d \tau}
\\ \nonumber
&\geq -\frac{2\alpha^{20}}{C} \Big( A^1_{t+N+n_j,t+N+n_{j+1}}(\gamma_1)+\frac{b^2 (n_{j+1}-n_j)}{2} \Big)
\\ \nonumber
&- \frac{C}{2\alpha^{20}} \Bigg( \int\limits_{t+N+n_j}^{t+N+n_j+\alpha^{40}}{ \Big(\frac{\partial G}{\partial x}(\gamma_1(\tau),\tau)-\frac{\partial G}{\partial x}(\gamma_1(\tau),t+N+n_{j+1}) \Big)^2 d \tau}
\\ \nonumber
&+ \int\limits_{t+N+n_j+\alpha^{40}}^{t+N+n_{j+1}}{ \Big(\frac{\partial G}{\partial x}(\gamma_1(\tau),\tau)-\frac{\partial G}{\partial x}(\gamma_1(\tau),t+N+n_{j+1}) \Big)^2 d \tau} \Bigg).
\end{eqnarray}
The first term of the right-hand side can be estimated by observing that the restriction of $\gamma_1$ to $[t+N+n_j,\ t+N+n_{j+1}]$ is still a minimizer, and that $A^3_{t+N+n_j,t+N+n_{j+1}}$, which only depends on the endpoint of the curve at $t+N+n_j$, is the same for $\gamma_1$ and $\gamma_3$.
\\ \indent
On the other hand, the second and the third terms of the right-hand side can be estimated by using (\ref{N2Nwhite}). Thus we obtain that
\begin{eqnarray} \nonumber
R &\geq -\frac{2\alpha^{20}}{C} \Big( - R+ A^1_{t+N+n_j,t+N+n_{j+1}}(\gamma_3) + A^2_{t+N+n_j,t+N+n_{j+1}}(\gamma_3) 
\\ \nonumber
&+\frac{b^2 N}{2} \Big) - \frac{C}{2\alpha^{20}} \Big( \frac{\alpha^{40} C^2}{4}+(n_{j+1}-n_j-\alpha^{40}) \alpha^{80} \Big)
\\ \nonumber
&\geq -\frac{2\alpha^{20}}{C} \Big( - R+\frac{1}{8(n_{j+1}-n_j)}
\\ \nonumber
&+ \int\limits_{t+N+n_j}^{t+N+n_{j+1}}{\dot{\gamma_3}(\tau) \Big(\frac{\partial G}{\partial x}(\gamma_3(\tau),\tau)-\frac{\partial G}{\partial x}(\gamma_3(\tau),t+N+n_{j+1})\Big) d \tau}
\\ \nonumber
&+\frac{b^2 N}{2} \Big) - \alpha^{20} C^3.
\end{eqnarray}
Consequently,
\begin{eqnarray}  \nonumber
R &\geq  \frac{2\alpha^{20}}{C}R-\frac{\alpha^{20}}{4C} 
\\ \nonumber
& -\frac{2\alpha^{20}}{C} \Big( b+\frac{1}{2} \Big) \int\limits_{t+N+n_j}^{t+N+n_{j+1}} { \Big|\frac{\partial G}{\partial x}(\gamma_3(\tau),\tau) }- \frac{\partial G}{\partial x}(\gamma_3(\tau),t+N+n_{j+1})\Big| d \tau 
\\ \nonumber
&- \frac{b^2 N \alpha^{20}}{C} - \alpha^{20} C^3.
\end{eqnarray}
Using (\ref{C}), (\ref{alphabis}), (\ref{NN}), and (\ref{N2Nwhite}), we get
\begin{eqnarray} \nonumber
R &\geq \frac{2\alpha^{20}}{C}R -\frac{\alpha^{20}}{4C} - \Big( b+\frac{1}{2} \Big)N \alpha^{20}-\frac{b^2 N \alpha^{20}}{C} -\alpha^{20}C^3
\\ \nonumber
&\geq \frac{2\alpha^{20}}{C}R - N \alpha^{20} (b+1)^2 - \Big( C^3+\frac{1}{4C}\Big) \alpha^{20}
\\ \nonumber
&\geq \frac{2\alpha^{20}}{C}R - \frac{4 \alpha^9}{10}-\frac{2\alpha^{17}}{10^3} \geq \frac{2\alpha^{20}}{C}R-\frac{\alpha^9}{2}.
\end{eqnarray}
Consequently,
\begin{eqnarray} \nonumber
&R \geq - \Big( 1-\frac{2\alpha^{20}}{C} \Big)^{-1} \frac{\alpha^9}{2} \geq -\alpha^{
9}.
\end{eqnarray}
By (\ref{nj}), it follows that for $j \in [1,N'-2]$ we have
\begin{eqnarray} \nonumber
&A_{t+N+n_j,t+N+n_{j+1}}(\gamma_2)-A_{t+N+n_j,t+N+n_{j+1}}(\gamma_1)
\\ \nonumber
&= (A^1_{t+N+n_j,t+N+n_{j+1}}(\gamma_2)-A^1_{t+N+n_j,t+N+n_{j+1}}(\gamma_1))
\\ \nonumber
&+ (A^2_{t+N+n_j,t+N+n_{j+1}}(\gamma_2)-A^2_{t+N+n_j,t+N+n_{j+1}}(\gamma_1))
\\ \nonumber
&+ (A^3_{t+N+n_j,t+N+n_{j+1}}(\gamma_2)-A^3_{t+N+n_j,t+N+n_{j+1}}(\gamma_1))
\\ \nonumber
&\leq (0-0)+\Big(b(C \alpha^{40}/2+N\alpha^{40})-(-\alpha^{
9})\Big)
\\ \label{kN12N2}
&+\Big((m_1+\alpha^6/10+\alpha^{40})-(m_1+\alpha^2-\alpha^{40}) \Big) \leq -\frac{\alpha^2}{2}.
\end{eqnarray}
Here, the estimate of $A^2_{t+N+n_j,t+N+n_{j+1}}(\gamma_2)$ follows from (\ref{N2Nwhite}).
\\ \indent
Similarly, since $A^3_{t+N,t+N+n_1}(\gamma_2)=A^3_{t+N,t+N+n_1}(\gamma_1)$, we have
\begin{eqnarray} \nonumber
A_{t+N,t+N+n_1}(\gamma_2)-A_{t+N,t+N+n_1}(\gamma_1)
\\ \label{kN}
\leq \frac{1}{8} +A^2_{t+N,t+N+n_1}(\gamma_2)+\alpha^{9} \leq 1,
\end{eqnarray}
and
\begin{eqnarray} \nonumber
&A_{t+N+n_{N'-1},t+2N}(\gamma_2)-A_{t+N+n_{N'-1},t+2N}(\gamma_1)
\\ \label{k2N1}
&\leq \frac{1}{8}+A^2_{t+N+n_{N'-1},t+2N}(\gamma_2)+\alpha^{9} + C \leq 2C.
\end{eqnarray}
Here, we get
\begin{eqnarray} \nonumber
&A^2_{t+N,t+N+n_1}(\gamma_2),A^2_{t+N+n_{N'-1},t+2N}(\gamma_2) \leq (b+1/2)(C \alpha^{40}/2+N\alpha^{40})
\end{eqnarray}
in the same way as for the estimate of $A^2_{t+N+n_j,t+N+n_{j+1}}(\gamma_2)$ above.
\\ \indent
It remains to add together the inequalities (\ref{kN12N2}-\ref{k2N1}). Using (\ref{alphabis}) and  (\ref{NN}) we get
\begin{eqnarray} \nonumber
&A_{t+N,t+2N}(\gamma_2)-A_{t+N,t+2N}(\gamma_1)
\\ \nonumber
&\leq 2C+1-(N'-2)\frac{\alpha^2}{2} \leq 2C+1-\frac{1}{2\alpha} < 0.
\end{eqnarray}
This inequality is in contradiction with the fact that $\gamma_1$ is a minimizer. This proves the lemma's assertion. $\square$

\begin{lemm} \label{white2}
For $\omega \in E_1 \cap E_2$, if for some minimizer $\gamma_1$ on $[t,t+N+n_j],\ 1 \leq j \leq N'-1$, $y=\gamma_1(t+N+n_j)$ satisfies:
\begin{eqnarray} \nonumber
&-\tilde{F}_1(y) \leq m_1+\alpha^2,
\end{eqnarray}
then we have
\begin{eqnarray} \nonumber
&\gamma_1(t+N) \in I_1(\alpha).
\end{eqnarray}
\end{lemm}

\textbf{Proof:} In the same way as in the proof of Lemma \ref{kick2}, we consider a \enquote{bad} minimizer $\gamma_1$. Without loss of generality, we assume that
\begin{equation} \label{assumwhite2}
-\tilde{F}_1[\gamma_1(t+N+n_{j'})] > m_1+\alpha^2,\ 1 \leq j' < j.
\end{equation}
We define $\gamma_2$ with the same endpoints as $\gamma_1$ in the same way as in the proof of Lemma \ref{kick2}, i.e. such that $\dot{\gamma_2}^b=0$ between $t+N$ and $t+N+n_j$, linear between $t$ and $t+N$, and satisfying $|\dot{\gamma_2}^b| \leq \frac{1}{2N}$. We get
\begin{eqnarray} \nonumber
&A^3_{t,t+N}(\gamma_2)=A^3_{t,t+N}(\gamma_1).
\\ \nonumber
&A^1_{t,t+N}(\gamma_2)-A^1_{t,t+N}(\gamma_1) \leq \frac{N}{8N^2}-0 \leq \frac{\alpha^{10}}{16}.
\\ \nonumber
&A^2_{t,t+N}(\gamma_2)\leq \Big(b+\frac{1}{2N} \Big) \int\limits_{t}^{t+N}{\Big| \frac{\partial G}{\partial x}(\gamma_2(\tau),\tau)-\frac{\partial G}{\partial x}(\gamma_2(\tau),t+N) \Big| d \tau} \leq \alpha^{29}.
\end{eqnarray}
The last inequality follows from (\ref{alphabis}), (\ref{NN}), and (\ref{0Nwhite}).
\\ \indent
To estimate the quantity
\begin{eqnarray} \nonumber
&R=A^2_{t,t+N}(\gamma_1),
\end{eqnarray}
we proceed in the same way as for $A^2_{t+N+n_j,t+N+n_{j+1}}(\gamma_1)$ in Lemma \ref{white1}. Namely, we consider a straight line $\gamma_3$ with the same endpoints as $\gamma_1|_{[t,t+N]}$ satisfying $|\dot{\gamma_3}^b| \leq 1/2N$. We have
\begin{eqnarray} \nonumber
R &\geq -2\alpha^{20} \int\limits_{t}^{t+N}{\Big( \frac{\dot{\gamma_1}(\tau)}{2}\Big)^2 d \tau}
\\ \nonumber
&- \frac{1}{2\alpha^{20}} \int\limits_{t}^{t+N}{\Big(\frac{\partial G}{\partial x}(\gamma_1(\tau),\tau)-\frac{\partial G}{\partial x}(\gamma_1(\tau),t+N)\Big)^2 d\tau}
\\ \nonumber
&\geq -2\alpha^{20} \int\limits_{t}^{t+N}{ \frac{(\dot{\gamma_1^b}(\tau))^2+b^2}{2} d \tau}
\\ \nonumber
&- \frac{1}{2\alpha^{20}} \int\limits_{t}^{t+N}{\Big(\frac{\partial G}{\partial x}(\gamma_1(\tau),\tau)-\frac{\partial G}{\partial x}(\gamma_1(\tau),t+N)\Big)^2 d\tau}.
\end{eqnarray}
Since a restriction of $\gamma_1$ is still a minimizer, we get
\begin{eqnarray} \nonumber
R &\geq -2\alpha^{20} \Big( A^1_{t,t+N}(\gamma_3)+A^2_{t,t+N}(\gamma_3)-R+\frac{b^2 N}{2} \Big)
\\ \nonumber
&- \frac{1}{2\alpha^{20}} \int\limits_{t}^{t+N}{\Big(\frac{\partial G}{\partial x}(\gamma_1(\tau),\tau)-\frac{\partial G}{\partial x}(\gamma_1(\tau),t+N)\Big)^2 d\tau}
\\ \nonumber
&\geq -2\alpha^{20} \Big( \frac{N}{8N^2} + \Big(b+\frac{1}{2N} \Big)
\\ \nonumber
&\times \int\limits_{t}^{t+N}{\Big|\frac{\partial G}{\partial x}(\gamma_3(\tau),\tau)-\frac{\partial G}{\partial x}(\gamma_3(\tau),t+N) \Big| d \tau}
\\ \nonumber
&-R + \frac{b^2 N}{2} \Big) - \frac{N}{2\alpha^{20}} \alpha^{80}
\\ \nonumber
&\geq  2\alpha^{20} R-2 \alpha^{20} \Big( \frac{\alpha^{10}}{16}+(b+1)N\alpha^{40}+\frac{2b^2}{\alpha^{10}}\Big) -\frac{N\alpha^{60}}{2}
\\ \nonumber
&\geq 2\alpha^{20} R-(5b^2+1) \alpha^{10} \geq 2\alpha^{20} R-\frac{\alpha^9}{2}.
\end{eqnarray}
Therefore
\begin{eqnarray} \nonumber
&A^2_{t,t+N}(\gamma_1) \geq -\alpha^{9}.
\end{eqnarray}
On the other hand, we have
\begin{eqnarray} \nonumber
&A^1_{t+N,t+N+n_j}(\gamma_2)-A^1_{t+N,t+N+n_j}(\gamma_1) \leq 0.
\end{eqnarray}
By definition, the action difference
\begin{eqnarray} \nonumber
&U=A_{t,t+N+n_j}(\gamma_2)-A_{t,t+N+n_j}(\gamma_1)
\end{eqnarray}
satisfies
\begin{eqnarray} \nonumber
U&= (A^1_{t,t+N}(\gamma_2)-A^1_{t,t+N}(\gamma_1))+(A^1_{t+N,t+N+n_j}(\gamma_2)-A^1_{t+N,t+N+n_j}(\gamma_1))
\\ \nonumber
&+ (A^2_{t,t+N}(\gamma_2)-A^2_{t,t+N}(\gamma_1))+ (A^3_{t,t+N}(\gamma_2)-A^3_{t,t+N}(\gamma_1))
\\ \nonumber
&+ (A^2_{t+N,t+N+n_j}(\gamma_2)+A^3_{t+N,t+N+n_j}(\gamma_2)
\\ \nonumber
&-A^2_{t+N,t+N+n_j}(\gamma_1)-A^3_{t+N,t+N+n_j}(\gamma_1)).
\end{eqnarray}
Consequently,
\begin{eqnarray}  \nonumber
U &\leq \frac{\alpha^{10}}{16}+0+\Big( \alpha^{29}+\alpha^{9} \Big)+0+ (A^2_{t+N,t+N+n_j}(\gamma_2)+A^3_{t+N,t+N+n_j}(\gamma_2)
\\ \nonumber
&-A^2_{t+N,t+N+n_j}(\gamma_1)-A^3_{t+N,t+N+n_j}(\gamma_1))
\\ \nonumber
&\leq  2\alpha^{9}+(A^2_{t+N,t+N+n_j}(\gamma_2)-A^2_{t+N,t+N+n_j}(\gamma_1)
\\ \label{quasifinal}
&+A^3_{t+N,t+N+n_j}(\gamma_2)-A^3_{t+N,t+N+n_j}(\gamma_1)].
\end{eqnarray}
In the same way as previously, we get
\begin{eqnarray} \nonumber
&A^2_{t+N,t+N+n_j}(\gamma_2) \leq bj \Big( \frac{C\alpha^{40}}{2}+N\alpha^{40} \Big) \leq 5bN' \alpha^{30} \leq \alpha^{26}.
\end{eqnarray}
The estimates of $A^2_{t+N+n_{j'},t+N+n_{j'+1}}(\gamma_1),\ 0 \leq j' < j$ in Lemma \ref{white1} still hold in our case. Therefore
\begin{eqnarray} \nonumber
&A^2_{t+N,t+N+n_j}(\gamma_1) \geq -N' \alpha^9 \geq -2 \alpha^6.
\end{eqnarray}
By (\ref{nj}) and (\ref{assumwhite2}), for $1 \leq j' \leq j-1$ we get
\begin{eqnarray} \nonumber
&A^3_{t+N+n_{j'},t+N+n_{j'+1}}(\gamma_2)-A^3_{t+N+n_{j'},t+N+n_{j'+1}}(\gamma_1)
\\ \nonumber
&\leq (m_1+\alpha^2+\frac{\alpha^6}{10}+\alpha^{40}) - (m_1+\alpha^2-\alpha^{40})\leq \alpha^6.
\end{eqnarray}
Finally, since we have supposed that $\gamma_1(t+N) \notin I_1(\alpha)$, we have
\begin{eqnarray} \nonumber
&A^3_{t+N,t+N+n_1}(\gamma_2)-A^3_{t+N,t+N+n_1}(\gamma_1)
\\ \nonumber
&\leq (m_1+\alpha^2+\frac{\alpha^6}{10}+\alpha^{40}) - (m_1+\alpha-\alpha^{40}) \leq -\alpha/2.
\end{eqnarray}
Combining all these inequalities with (\ref{quasifinal}) we get
\begin{eqnarray} \nonumber
&U \leq 2\alpha^{9}+\alpha^{26}+2\alpha^{6}+(N'-1)\alpha^{6} -\alpha/2 < 0.
\end{eqnarray}
We have a contradiction with the fact that $\gamma_1$ is a minimizer. This proves the lemma's assertion. $\square$

\ack{}

We would like to thank R. Iturriaga and S. Kuksin for numerous discussions.
A big part of the work was carried out while both authors were visiting Observatoire de la C{\^o}te d'Azur in Nice. We are very grateful to U. Frisch, J. Bec and other colleagues
from the Observatoire for their hospitality.

\bigskip


\bigskip

\begin{thebibliography}{9}

\bibitem{BK}
J.~Bec, K.~Khanin, \textit{Burgers Turbulence}, Phys. Rep. 447(1-2), 2007, 1-66.

\bibitem{EKMS}
Weinan E, K.~Khanin, A.~Mazel, Ya.~Sinai, 
\textit{Invariant measures for Burgers equation with stochastic forcing}, Ann. of Math. 101(3), 2000, 877-960. 

\bibitem{IK}
R.~Iturriaga, K.~Khanin, \textit{Burgers Turbulence and Random Lagrangian Systems}, Commun. Math. Phys. 232 (3), 2003, 377-428.
%

\bibitem{Rud}
W.~Rudin, \textit{Real and Complex Analysis}, McGraw-Hill, 1987.

\end{thebibliography}
\end{document}